\documentclass[10pt]{article}
\usepackage{amssymb}
\usepackage{graphicx}

\def\proof{\noindent{\bf Proof.\quad}}
\def\qed{\hspace*{\fill}\vrule height6pt width4pt depth0pt\medskip}

\newtheorem{lemma}{Lemma}
\newtheorem{thm}[lemma]{Theorem}
\newtheorem{pro}[lemma]{Proposition}

\begin{document}
\baselineskip = 10pt
%

\title{Edge-colorings of $K_{m,n}$ which Forbid Multicolored Cycles
\thanks {Research supported in part by NSC 97-2115-M-009-011-MY3.}}

\author{Hung-Lin Fu, Yuan-Hsun Lo and Ryo-Yu Pei\\ \\
{\small \it Department of Applied Mathematics}\\
{\small \it National Chiao Tung University}\\
{\small \it 1001 Ta Hsueh Road, Hsinchu, Taiwan 300, ROC}}
\date{}
\maketitle

\begin{abstract}
A subgraph in an {\it edge-colored} graph is {\it multicolored} if all its edges receive distinct colors. In this paper, we study the proper edge-colorings of the complete bipartite graph $K_{m,n}$ which forbid multicolored cycles. Mainly, we prove that (1) for any integer $k\geq 2$, if $n\geq 5k-6$, then any properly $n$-edge-colored $K_{k,n}$ contains a multicolored $C_{2k}$, and (2) determine the order of the properly edge-colored complete bipartite graphs which forbid multicolored $C_6$.
\end{abstract}

\section{Introduction}

Throughout this paper, all terminologies and notations on graph theory can be referred to the textbook by D. B. West \cite{W}. A ({\it proper}) {\it k-edge-coloring} of a graph $G$ is a mapping from $E(G)$ into a set of colors $\{1,2,\ldots,k\}$ such that incident edges of $G$ receive distinct colors. It's well-known that $K_{m,n}$ has an $n$-edge-coloring where $m$ and $n$ are natural number and $n\geq m$. A subgraph in an edge-colored graph is said to be $multicolored$ (or $rainbow$) if no two edges have the same color.

Suppose $G$ is a finite simple graph and $H$ is a subgraph of $G$. An edge-coloring of $G$ $forbids$ multicolored $H's$ (copies of $H$) if each copy of $H$ in $G$ has two edges with the same color. A. Gouge et al. \cite{GHJNP} discussed the edge-colorings of $K_n$ that forbid multicolored $K_3's$ and thus all multicolored cycles. In this paper, motivated by their works, we consider the edge-colorings of $K_{m,n}$, $n\geq m$, which forbid multicolored (even) cycles. Actually, given an integer $k$, we want to know for what natural numbers $n$ and $m$, there always exists a multicolored $C_{2k}$ somewhere in any $n$-edge-colored $K_{m,n}$. For $k\geq 2$, we define the {\it forbidding multicolored 2k-cycles set}, $FMC(2k)$ in short, by the ordered pair $(m,n)\in FMC(2k)$ if there exists an $n$-edge-coloring of $K_{m,n}$ that forbids multicolored $2k$-cycles. Since $m<k$ or $n<2k$ gives trivial results, we only consider $m\geq k$ and $n\geq 2k$ in the set $FMC(2k)$.

Firstly, it is impossible to forbid multicolored $4$-cycles in any $n$-edge-coloring of $K_{m,n}$ where $2\leq m\leq n$ and $n\geq 4$.

\begin{pro}\label{C4} $FMC(4)=\phi$.
\end{pro}
\proof It suffices to show that there exists a multicolored $C_4$ in any properly edge colored $K_{2,4}$. Let $\varphi$ be a properly edge coloring of $K_{2,4}$ and $\{ u_1,u_2 \}$, $\{ v_1,v_2,v_3,v_4 \}$ be the two partite sets. For convenience, let $C=\{1,2,\ldots \}$ be the color set. Without loss of generality, assume $\varphi(u_1v_1)=1$ and $\varphi(u_2v_1)=2$. There must be one vertex $v_i$, where $2\leq i\leq 4$, such that $\varphi(u_1v_i),\varphi(u_2v_i)\notin \{ 1,2\}$. Thus $u_1-v_1-u_2-v_i-u_1$ is the desired multicolored $C_4$. \qed

\section{Forbidding Multicolored $2k$-cycles}

By Proposition \ref{C4}.

Let $S$ be an $n$-set. A {\it latin rectangle} of order $m\times n$, $m\leq n$, based on $S$ is an $m\times n$ array in which every element of $S$ is arranged such that each one occurs at most once in each row and each column. A {\it latin square} of order $n$ based on $S$ is a latin rectangle of order $n\times n$. A {\it partial latin square} of order $r$, $r<n$, based on $S$ is an $r\times r$ array in which every element of $S$ is arranged such that each one occurs at most once in each row and each column. In this paper, we use $\mathbb{Z}_n=\{0,1,2,\ldots,n-1\}$ for the $n$-set $S$. For example,
\begin{tabular}{|c|c|c|}
\hline 0 & 1 & 2 \\
\hline 2 & 0 & 1 \\
\hline
\end{tabular}
is a latin rectangle of order $2\times 3$ based on $\mathbb{Z}_3$;
\begin{tabular}{|c|c|}
\hline 0 & 1 \\
\hline 1 & 0 \\
\hline
\end{tabular}
is a latin square of order 2 based on $\mathbb{Z}_2$; and
\begin{tabular}{|c|c|}
\hline 0 & 1 \\
\hline 2 & 0 \\
\hline
\end{tabular}
is a partial latin square of order 2 based on $\mathbb{Z}_3$. In particular, the size of a partial latin square $L$, denoted by $|L|$, is the number of elements of $S$ actually appearing in $L$.

\begin{figure}[h]
    \begin{center}
        \includegraphics[scale=0.6]{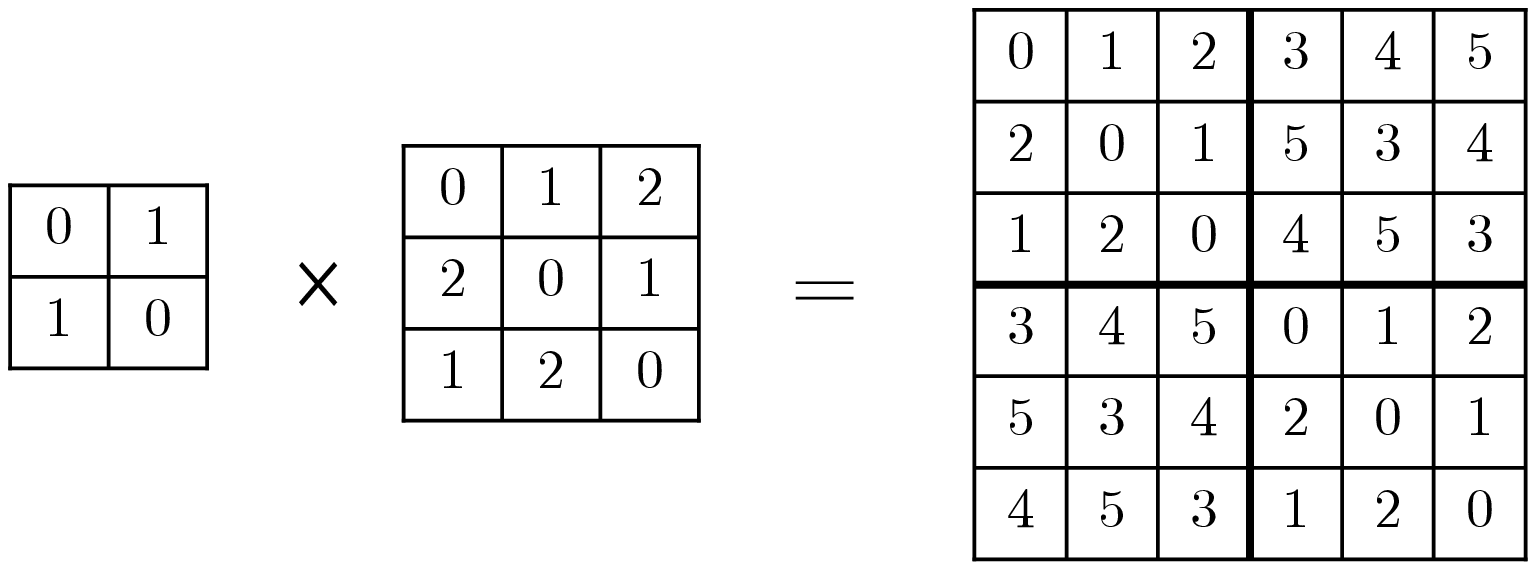}
    \end{center}
    \caption{\label{product} The direct product of $L$ and $M$}
\end{figure}

For convenience, a latin square of order $n$ based on $\mathbb{Z}_n$ is denoted $L=[~l_{i,j}~]$ where $l_{i,j}\in \mathbb{Z}_n$ and $i,j \in \mathbb{Z}_n$. Let $L=[~l_{i,j}~]$ and $M=[~m_{i,j}~]$ be two latin squares of order $s$ and $t$ respectively. Then the direct product of $L$ and $M$ is a latin square of order $s\cdot t$ : $L\times M=[~h_{i,j}~]$, where $h_{x,y}=t\cdot l_{a,b}+m_{c,d}$ provided that $x=ta+c$ and $y=tb+d$. For instance, let $L$ and $M$ be two latin square of order $2$ and $3$ respectively, then $L\times M$ is a latin square of order 6 based on $\mathbb{Z}_6$, as in Figure \ref{product}.

A {\it transversal} of a latin square of order $n$ is a set of $n$ cells with exactly one in each
row and each column and containing exactly $n$ elements. For example, in $M$ (Figure
\ref{product}), $\{m_{0,0},m_{1,2},m_{2,1}\}$ is a transversal. Similarly, a transversal of a
partial latin square of order $r$ based on an $n$-set is set of $r$ cells with exactly one in each
row and each column and containing exactly $r$ elements. It's easy to check there is no transversal
in $L\times M$ (Figure \ref{product}). For more information on latin squares, we refer to
\cite{DK}.

Let $L=[l_{i,j}]$ be an $m\times n$ latin rectangle. There is a corresponding relationship between $L$ and an $n$-edge-colored $K_{m,n}$. Let $\{ u_0,u_1,\ldots, u_{m-1} \}$ and $\{ v_0,v_1,\ldots,v_{n-1} \}$ be two partite sets of $K_{m,n}$, and the edge $u_iv_j$ is colored with $l_{i,j}$ for each $0\leq i\leq m-1$, $0\leq j\leq n-1$, then we have an $n$-edge-colored $K_{m,n}$ and vice versa. Now, we have

\begin{thm}\label{C2k-lowerbound} If $k$ is odd, then $(m,2k)\in FMC(2k)$ for $k\leq m\leq 2k$.
\end{thm}

\begin{figure}[h]
    \begin{center}
        \includegraphics[scale=0.6]{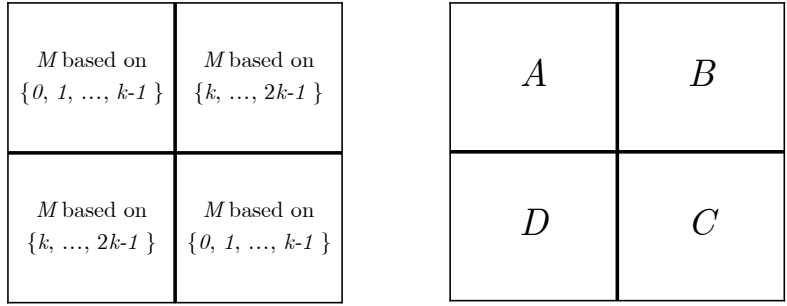}
    \end{center}
    \caption{\label{f2k} $L_2 \times M$ and the four copies of $M$}
\end{figure}

\proof It suffices to find a $2k$-edge-coloring of $K_{2k,2k}$ which forbids multicolored $C_{2k}$. Let $L_2$ be the latin square of order 2 in Figure \ref{product} and $M$ be a latin square of order $k$. Notice that $L_2\times M$ is formed by four latin squares of order $k$, two of them based on $\mathbb{Z}_k$ and other two based on $\mathbb{Z}_{2k} \setminus \mathbb{Z}_k$. For convenience, name the four squares $A,B,C$ and $D$ clockwise from the top-left one, see Figure \ref{f2k}.

Let $\varphi$ be the $2k$-edge-coloring of $K_{2k,2k}$ obtained by $L_2\times M$. Suppose it contains a multicolored $C_{2k}$. Let $a,b,c,$ and $d$ be the numbers of cells in $A,B,C,$ and $D$, respectively, corresponding to the edges of the multicolored cycle. Then $a+b$ is a sum of the degrees, on the cycle, of some of the vertices on the cycle, so $a+b$ is even. Similarly, $b+c$ is even. Therefore, $a+c$ is even. But since all $2k$ colors $0,1,\ldots,2k-1$ must appear on the edges of the cycle, $a+c=k$, odd. This contradiction completes the proof. \qed

The following result provides an upper bound of the order of complete bipartite graphs to forbid multicolored $2k$-cycles.

\begin{thm}\label{C2k-upperbound} For any integer $k\geq 2$, if $n\geq 5k-6$, then any $n$-edge-colored $K_{k,n}$ contains a multicolored $C_{2k}$.
\end{thm}

\proof Let $\varphi$ be an $n$-edge-coloring of $K_{k,n}$ and the partite sets be $A=\{ a_1,a_2,\ldots, a_k \}$ and $B=\{ b_1,b_2,\ldots, b_n \}$. Let $P=a_1b_1a_2\cdots b_{t-1}a_t$ be the longest multicolored path whose endpoints lie on $A$. Suppose $t<k$. Assume $C$ is the set of colors which appear on $P$. Note that $|C|=2t-2$. For each $i=1,\ldots,k$, define $S_i\subset B$ by $b\in S_i$ if $\varphi(a_ib)\in C$. Observe that $|S_t\cup S_{t+1}\cup \{b_1,b_2,\ldots,b_{t-1}\}|\leq 2(2t-2)+(t-1)-1=5t-6<5k-6\leq n$. Therefore, there exists a vertex $b\in \{b_t,b_{t+1},\ldots,b_n\}$ such that $\varphi(a_tb), \varphi(a_{t+1}b)\notin C$, a contradiction. Then $t\geq k$. By the fact that a longest path in $K_{k,n}$ with end vertices in $A$ is of length $2k-2$, we have $t=k$.

We have that $|S_1|,|S_k|\leq 2k-2$ and $b_1\in S_1,b_{k-1}\in S_k$. Hence, $|S_1\cup S_k\cup \{b_1,\ldots,b_{k-1}\}|\leq 5k-7$. Since $n\geq 5k-6$, there exists a vertex $b\in B$ such that $\varphi(a_1b),\varphi(a_kb)\notin C$. Therefore, a multicolored $C_{2k}$ is found. \qed

\section{Determining $FMC$(6)}

By Theorem \ref{C2k-upperbound}, if $(m,n)\in FMC(6)$, then we have $3\leq m\leq n$ and $n=6,7,8$. The case $n=6$ was done in Theorem \ref{C2k-lowerbound}, so we consider $n=7$ and $8$ in the following.

Let $L$ be the corresponding latin rectangle of an $n$-edge-colored $K_{m,n}$. If there is a multicolored $C_6$ somewhere, then there exists a $3\times 3$ partial latin square which contains two disjoint transversals using exactly $6$ symbols in $L$.

\begin{pro}\label{7elements} Let $L$ be a partial latin square of order 3 with $|L|=7$. Then, there is no multicolored $C_6$ in its corresponding $K_{3,3}$ if and only if it contains a latin subsquare of order $2$.
\end{pro}

\proof It suffices to consider the necessity since the sufficiency is clearly true. Suppose $L$ contains no latin subsquares of order $2$. If there is one element appearing $3$ times, then the other $6$ elements form a multicolored $C_6$. Therefore, assume that there are two elements, say $1,2$, appearing twice respectively. Without loss of generality, let the two $1$'s be arranged at the diagonal in the first two rows. Then $2$ occurs in the third column or the third row. Omitting this cell and one of the cells labeled $1$, the 6 of the the remaining cells will provide a multicolored $C_6$, a contradiction. \qed

\begin{pro}\label{6elements} Let $L$ be a partial latin square of order 3 with $|L|=6$. There does not exist a multicolored $C_6$ in its corresponding $K_{3,3}$ if one of the following conditions occurs:
\begin{enumerate}
\item There exist $2$ columns (or rows) in $L$ using exactly $3$ elements.
\item Some element appears three times in $L$.
\item $L$ contains a latin subsquare of order $2$.
\end{enumerate}
\end{pro}

\proof Since there are only 6 elements, if there exists a multicolored $C_6$, all elements should appear in the two disjoint transversals. In case 1, the elements of the third column (or row) can not all appear. In case 2, that element can not appear only once in any two disjoint transversals. In case 3, the argument is similar to the proof of Proposition \ref{7elements}. \qed

\begin{lemma}\label{K-m8} For $3\leq m\leq 8$, $(m,8)\in FMC(6)$.
\end{lemma}

\begin{figure}[h]
    \begin{center}
        \includegraphics[scale=0.6]{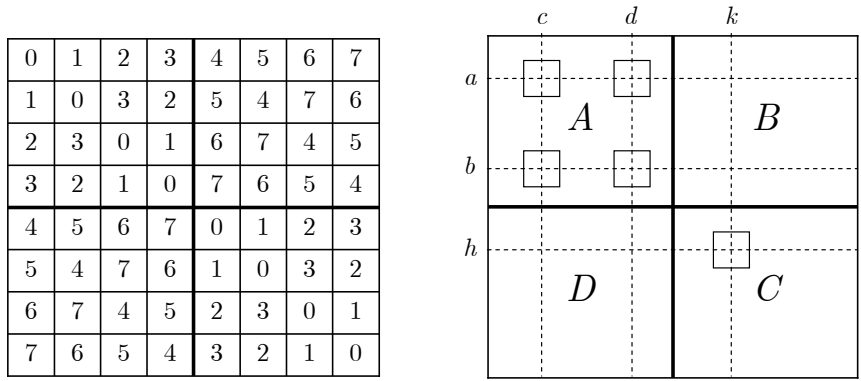}
    \end{center}
    \caption{\label{88} $L_2\times L_2\times L_2$}
\end{figure}

\proof It suffices to prove the claim for $m=8$. Let $L_2$ be the latin sqaure of order 2 in Figure \ref{product}. Let $L=L_2\times L_2\times L_2$, a latin square of order $8$ based on $\mathbb{Z}_8$. For convenience, name the four copies $A,B,C$ and $D$ of $L_2\times L_2$ as in Figure \ref{88}.

Suppose that there are $6$ cells whose entries induce a multicolored $C_6$. Let $L'$ be the $3\times 3$ partial latin square which contains the $6$ cells. It is easy to see that any $2\times 3$ partial latin rectangle in $L_2\times L_2$ ($A$ or $B$ or $C$ or $D$) contains a latin subsquare of order $2$. By Proposition \ref{7elements}, we can assume that $L'$ traverses all four copies of $L_2\times L_2$. Without loss of generality, suppose there are $4$ cells of $L'$ located in $A$. Let the $4$ cells be $(a,c),(a,d),(b,c),(b,d)$, and the only one cell located on $C$ be $(h,k)$, where $0\leq a,b,c,d\leq 3$ and $4\leq h,k\leq 7$ (Figure \ref{88}). By Proposition \ref{7elements} and Proposition \ref{6elements}, $l_{a,c}\neq l_{b,d}$ or $l_{a,d}\neq l_{b,c}$, and thus the four elements are distinct. Assume that $l_{h,k}=l_{a,c}$. This implies $l_{a,k}=l_{h,c}$. Then we have a copy of $L_2$, a contradiction. Similarly, if $l_{h,k}$ is any of the $l_{i,j}$, with $(i,j)$ being one of the 4 cells of $L'$ in $A$, then we have a contradiction. But $l_{h,k}$ must be one of these, since these 4 are distinct elements of $\{ 0,1,2,3\}$. \qed

\begin{lemma}\label{K-37} $(3,7)\in FMC(6)$. Furthermore, if $K_{3,7}$ is $7$-edge-colored such that it forbids multicolored $C_6$'s, there exists an induced $K_{3,3}$ using exactly $3$ colors.
\end{lemma}

\begin{figure}[h]
    \begin{center}
        \includegraphics[scale=0.7]{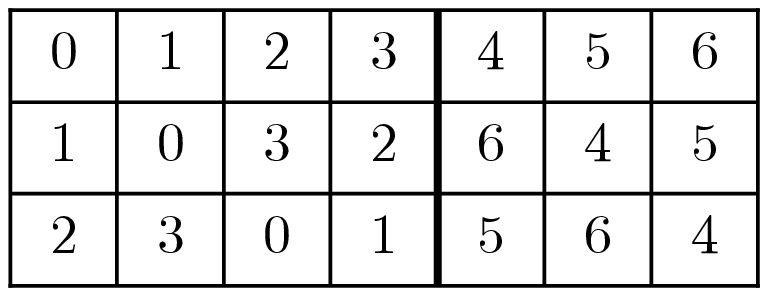}
    \end{center}
    \caption{\label{37} A $3\times 7$ latin rectangle}
\end{figure}

\proof Firstly, Figure \ref{37} gives a $3\times 7$ latin rectangle. It is not difficult to check its corresponding $7$-edge-coloring of $K_{3,7}$ induces no multicolored $C_6$ by Proposition \ref{7elements} and Proposition \ref{6elements}.

Secondly, given a $7$-edge-coloring of $K_{3,7}$ which forbids multicolored $6$-cycles and let $L$ be its corresponding latin rectangle. It suffices to show that $L$ contains a latin subsquare of order $3$. For convenience, let $C^i$ denote the set of elements in the $i$th column of $L$ where $i\in \mathbb{Z}_7$.

{\bf Claim.} There exist $i,j$ such that $C^i\cap C^j=\phi$.

Suppose for any $i\neq j$, $C^i\cap C^j\neq \phi$. Since each element occurs three times, we have $|C^i\cap C^j|=1$ for all $i\neq j$ under this assertion. Without loss of generality, let $C^0=\{0,1,2\}$ and $C^1=\{0,3,4\}$. Then 3 and 4 will each occur twice in the remaining five columns. So, there exists a $C^t$, where $2\leq t\leq 6$, such that $C^t\cap \{3,4\}=\phi$. This implies that the three columns $C^0,C^1$ and $C^t$ create a multicolored $C_6$ by Proposition \ref{7elements}, a contradiction.

\begin{figure}[h]
    \begin{center}
        \includegraphics[scale=0.7]{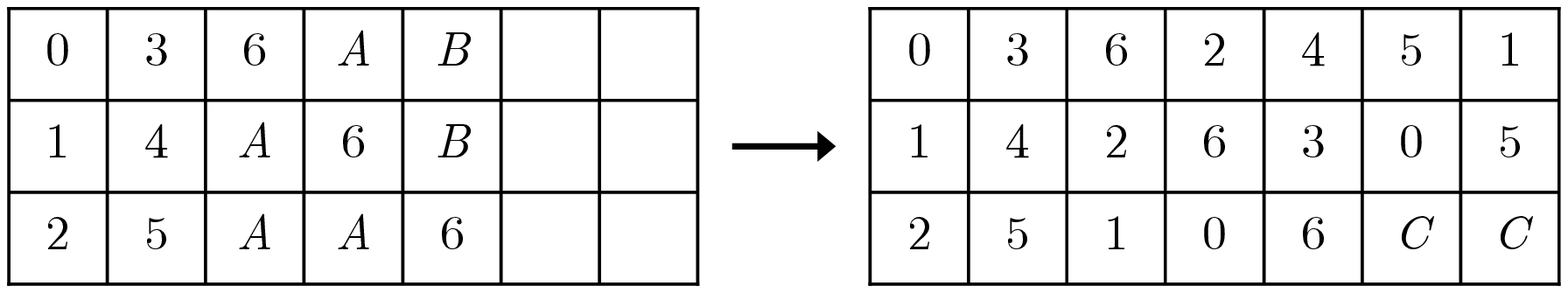}
    \end{center}
    \caption{\label{prove37} The $3\times 7$ latin rectangle}
\end{figure}

Thus, assume $C^0=\{ 0,1,2\}, C^1=\{ 3,4,5\}$ and $C^2,C^3,C^4$ contain the element 6. Note here that $|C^t\cap C^0|=2$ or $|C^t\cap C^1|=2$ for all $t=2,3,4$; otherwise, $C^0,C^1,C^t$ will create a multicolored $C_6$ by Proposition \ref{7elements}. Next, we want to claim $(C^2\cup C^3\cup C^4)\setminus \{6\}$ equals $C^0$ or $C^1$. Suppose the assertion is not true, without loss of generality, let $|C^2\cap C^0|=2, |C^3\cap C^0|=2$ and $|C^4\cap C^1|=2$. See the left rectangle in Figure \ref{prove37}: the elements in cell $A$ are from $\{ 0,1,2\}$ while the elements in cell $B$ are from $\{ 3,4,5\}$.

Proposition \ref{7elements} shows that the elements in the cells labelled $A$ and the cells labelled $B$ are uniquely determined; see the right hand side rectangle in Figure \ref{prove37}. Meanwhile, the elements in some cells of the last two columns are determined except cells denoted as $C$, which are filled with 3 and 4. No matter what the elements in $C$ are, $C^0,C^4$ and $C^5$ contain a multicolored $C_6$, a contradiction. Therefore, $(C^2\cup C^3\cup C^4)\setminus \{6\}$ equals $C^0$(or $C^1$). Hence, combining $C^5,C^6$ with $C^1$(or $C^0$), we have a latin square of order 3. \qed

Lemma \ref{K-37} will yield the following result.

\begin{pro}\label{K-m7} For any $7$-edge-coloring of $K_{m,7}$, $4\leq m\leq 7$, there exists a multicolored $C_6$.
\end{pro}

\proof

\begin{figure}[h]
    \begin{center}
        \includegraphics[scale=0.7]{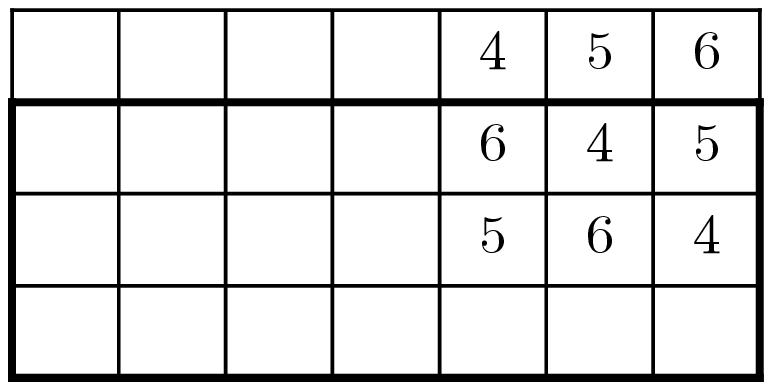}
    \end{center}
    \caption{\label{m7} The $4\times 7$ latin rectangle}
\end{figure}

It's sufficient to consider the case when $m=4$. Suppose there exists a 7-edge-colored $K_{4,7}$ which forbids multicolored $C_6$'s. Let $L$ be its corresponding $4\times 7$ latin rectangle. By Lemma \ref{K-37}, there exists a latin square of order 3 in the first three rows of $L$. Without loss of generality, we put the latin square of order 3 in the last three columns and let the symbols be $\{4,5,6\}$, see Figures \ref{m7}. Next, consider the last three rows. It's impossible to find another latin square of order 3. It contradicts Lemma \ref{K-37}. \qed

To sum up, we have the following conclusion.

\begin{thm}\label{FMC6} $FMC(6)=\{ (m,6)|~3\leq m\leq 6 \} \cup \{ (3,7) \} \cup \{ (m,8)|~3\leq m\leq 8 \}$.
\end{thm}

\newpage

\noindent{\bf \Large Acknowledgements}

The authors would like to express their gratitude to the referees for their careful reading and their many important comments that significantly improved the presentation of the article.

\end{document}